\documentclass{article}
\usepackage{amsmath}
\usepackage{amsfonts,amsthm,enumerate}%
\usepackage{amssymb}
\usepackage[pdftex]{graphicx}

\title{Logarithmically-small Minors and Topological Minors}
\author{Richard Montgomery\footnote{Department of Pure Mathematics and Mathematical Statistics, Centre for Mathematical Sciences, Wilberforce Road, Cambridge, CB3 0WB, UK. r.h.montgomery@dpmms.cam.ac.uk \newline \indent
Supported by the Engineering and Physical Sciences Research Council.}}

\newcounter{CTR}
\newtheorem{lemma}[CTR]{Lemma}
\newtheorem{corollary}[CTR]{Corollary}
\newtheorem{theorem}[CTR]{Theorem}
\newtheorem{prop}[CTR]{Proposition}
\newtheorem{defn}[CTR]{Definition}

\global\long\def\a{\alpha}
\global\long\def\b{\beta}
\global\long\def\c{\gamma}
\global\long\def\d{\delta}

\global\long\def\e{\varepsilon}

\global\long\def\PP{\mathcal{P}}
\global\long\def\QQ{\mathcal{Q}}

\global\long\def\de{\begin{defn}}
\global\long\def\fn{\end{defn}}
\global\long\def\cor{\begin{corollary}}
\global\long\def\ary{\end{corollary}}
\global\long\def\lem{\begin{lemma}}
\global\long\def\ma{\end{lemma}}
\global\long\def\arr{\begin{array}}
\global\long\def\ay{\end{array}}
\global\long\def\pr{\begin{proof}}
\global\long\def\oof{\end{proof}}

\begin{document}

\maketitle

\begin{abstract}
For every integer $t$ there is a smallest real number $c(t)$ such that any graph with average degree at least $c(t)$ must contain a $K_t$-minor (proved by Mader). Improving on results of Shapira and Sudakov, we prove the conjecture of Fiorini, Joret, Theis and Wood that any graph with $n$ vertices and average degree at least $c(t)+\e$ must contain a $K_t$-minor consisting of at most $C(\e,t)\log n$ vertices.

Mader also proved that for every integer $t$ there is a smallest real number $s(t)$ such that any graph with average degree larger than $s(t)$ must contain a $K_t$-topological minor. We prove that, for sufficiently large $t$, graphs with average degree at least $(1+\e)s(t)$ contain a $K_t$-topological minor consisting of at most $C(\e,t)\log n$ vertices. Finally, we show that, for sufficiently large $t$, graphs with average degree at least $(1+\e)c(t)$ contain either a $K_t$-minor consisting of at most $C(\e,t)$ vertices or a $K_t$-topological minor consisting of at most $C(\e,t)\log n$ vertices.
\end{abstract}

\section{Introduction}
A graph $H$ is a minor of a graph $G$ if a copy of $H$ can be formed from $G$ by contracting edges and deleting edges and vertices. This fundamental notion was introduced by Wagner~\cite{wag37b} in his study of Kuratowski's theorem characterising planar graphs. He proved that the four colour theorem is equivalent to the statement that every graph of chromatic number at least 5 contains $K_5$ as a minor~\cite{wag37}. Hadwiger subsequently conjectured that the statement holds with 5 replaced by any integer, and this remains one of the outstanding open questions in the subject. A deep conjecture of Wagner himself, that in any infinite collection of finite graphs there is one which is a minor of another, has, however, been resolved in a renowned series of papers by Robertson and Seymour~\cite{RS04}.

A principal result in the study of graph minors is that, for each integer $t$, a sufficiently large average degree, based only on $t$, forces a graph to contain the complete graph on $t$ vertices as a graph minor. That is, if $d(G)$ is the average degree of a graph $G$ and
\[
c(t)=\min\{c:d(G)\geq c \text{ implies that }G\text{ has a }K_t\text{-minor}\},
\]
then $c(t)$ is finite. Mader~\cite{Mad67,Mad68} showed that $c(t)$ does indeed take the minimum of the above set and that $c(t)\leq 2^{t-2}$, before improving the bound to $c(t)=O(t\log t)$. Kostochka~\cite{Kos84} and Thomason~\cite{AGT84} independently found the correct order, showing that $c(t)=\Theta(t\sqrt{\log t})$, where the lower bound comes from a disjoint union of dense random graphs. Thomason~\cite{AGT01} was later able to prove that $c(t)=(2\a+o(1))t\sqrt{\log t}$ for an explicit constant $\a$.

It is natural to ask how many vertices are needed to form the $K_t$-minor in graphs whose average degree exceeds $c(t)$. Fiorini, Joret, Theis and Wood~\cite{FJTW12} proved that an average degree of at least $2^{t-1}+\e$ in a graph with $n$ vertices suffices to force a $K_t$-minor constructed from $O_{\e,t}(\log n)$ vertices. That is to say they showed there was some constant $C(\e,t)$ depending on $\e$ and $t$ so that any graph with $n$ vertices and average degree at least $2^{t-1}+\e$ must contain a $K_t$-minor found by contracting at most $C(\e,t)\log n$ edges and deleting other edges and vertices. They conjectured that only an average degree of at least $c(t)+\e$ would be needed to force a $K_t$-minor consisting of $O_{\e,t}(\log n)$ vertices. As they remarked, the existence of graphs with constant average degree $d$ yet girth $C(d)\log n$ (either found by random methods~\cite{ES63} or explicit constructions~\cite{Cha03}) demonstrates that this would be tight up to the constant when $t\geq 3$. Shapira and Sudakov~\cite{SS12} used graph expansion to get near this conjecture, showing that each graph with average degree at least $c(t)+\e$ must contain a $K_t$-minor consisting of $O_{\e,t}(\log n\log\log n)$ vertices. Here we prove the full conjecture.

\begin{theorem} \label{logbound}
Let $\e>0$. If a graph $G$ with $n$ vertices has average degree at least $c(t)+\e$ then it contains a $K_t$-minor consisting of $O_{\e,t}(\log n)$ vertices.
\end{theorem}

A graph $G$ is a subdivision of a graph $H$ if it can be formed by replacing the edges of $H$ by disjoint paths. We say that $G$ contains an $H$-subdivision if an $H$-subdivision can be found by deleting vertices and edges from $G$. This is also known as $G$ containing $H$ as a topological minor. A well known result is that, for each integer $t$, a sufficiently large average degree, based only on $t$, forces a graph to contain $K_t$ as a subdivision. That is, if
\[
s(t)=\mbox{inf}\{c:d(G)\geq c \text{ implies that }G\text{ has a }K_t\text{-subdivision}\},
\]
then $s(t)$ is finite. More difficult than the previous case with minors, this was first proved by Mader~\cite{Mad67}, who later improved his result to $s(t)=O(2^t)$~\cite{Mad72}.

Koml\'os and Szemer\'edi~\cite{KSTI} used graph expansion in the sparse case and Szemer\'edi's Regularity Lemma in the dense case to prove that for each $\gamma>14$ there is some constant $c$ such that $s(t)\leq ct^2(\log t)^\gamma$. Bollob\'as and Thomason~\cite{BT98top} found the correct order, showing that $s(t)=\Theta(t^2)$ using the concept of linked graphs. Subsequently, Koml\'os and Szemer\'edi~\cite{KSTII} improved their previous methods to give an alternative proof. Currently it is known that
\[
(1+o(1))9t^2/64\leq s(t)\leq (1+o(1))10t^2/23,
\]
where the lower-bound is due to an example by \L uczak, and the upper-bound is due to a development of Koml\'os and Szemer\'edi's proof by K\"uhn and Osthus~\cite{KO06}.

As noted by Fiorini, Joret, Theis and Wood~\cite{FJTW12}, the methods used by Kostochka and Pyber in~\cite{KP88} can be adapted to show that, for every integer $t\geq 2$ and real $\e>0$, any graph on $n$ vertices with average degree at least $4^{t^2}+\e$ contains a $K_t$-subdivision consisting of $O_{\e,t}(\log n)$ vertices.

Developing our method for constructing small minors further, we will prove the following result for subdivisions.

\begin{theorem} \label{logboundtop} Let $\e>0$. For sufficiently large $t\geq t_0(\e)$ each graph with $n$ vertices and average degree at least $(1+\e)s(t)$ contains a $K_t$-subdivision consisting of $O_{\e,t}(\log n)$ vertices.
\end{theorem}

In Theorem \ref{logboundtop} we require slightly stronger conditions than the direct analogue of Theorem \ref{logbound}, which would see the requirement for $t$ to be sufficiently large removed and the degree bound reduced to $s(t)+\e$. These slightly stronger conditions reflect the additional difficulties in constructing a subdivision, but the direct analogue could plausibly hold.

Using very similar methods, we will prove the following theorem, which finds either a constant-sized minor or a logarithmic-sized subdivision in any graph with above the extremal edge density for graph minors.

\begin{theorem}\label{topsinminors} Let $\e>0$. For sufficiently large $t\geq t_0(\e)$ each graph with $n$ vertices and average degree at least $(1+\e)c(t)$ contains either a $K_t$-minor consisting of $O_{\e,t}(1)$ vertices or a $K_t$-subdivision consisting of $O_{\e,t}(\log n)$ vertices.
\end{theorem}

The main tool we will use to construct our small minors is graph expansion. A graph $G$ is typically said to be an expander if for every vertex set $A$ with size smaller than some $m$ the neighbourhood of $A$ has size at least $\lambda|A|$. Usually the values of $m$ and $\lambda$ are chosen to balance what conditions are needed to guarantee expansion against what properties can be derived from the expansion.

A central idea in the proof of Shapira and Sudakov~\cite{SS12} is that, roughly speaking, every graph contains an expander subgraph of almost the same average degree as the original graph (see Lemma \ref{key} for a version of this). This reduces the problem to that of finding small minors in expander graphs.

Our approach allows us to find small minors using a weaker expansion property in which only small to medium sized subsets expand. The minor is constructed by finding very many expanding subsets of moderate size, amongst which there must be some that can be used to build a minor. Our definition of expansion is thus weaker than in \cite{SS12} and the overall proof is hence shorter.

This approach also serves as the basis for constructing subdivisions, but further ideas and tools are needed. As the details are more intricate we defer an outline to Section \ref{sexp2}.

In Section \ref{sexp} we will set up some notation and outline the graph expansion concepts we will require, before proving Theorem \ref{logbound} in Section \ref{main}. In Section \ref{sexp2} we will strengthen the graph expansion concepts for small sets. In Section \ref{main2} we will describe units, how they can be found in graphs with expansion properties, and how they can be used to construct subdivisions. In Section \ref{main3} we will prove Theorems \ref{logboundtop} and \ref{topsinminors} using the results in Section \ref{main2}. In Section \ref{remarks} we will make some closing remarks.

In several lemmas about graphs of order $m$ we have omitted rounding symbols for functions of $m$ for neatness. This will not affect the proofs as the lemmas will only be applied for large $m$ (depending on $\e$ and $t$).

\section{Graph Expansion and Notation}\label{sexp}

\subsection{Basic Notation}

For a graph $G$ and a vertex subset $W\subset V(G)$, $G[W]$ is the subgraph of $G$ induced on the vertices of $W$. We will abbreviate $G[V(G)\setminus W]$ to $G-W$. We will denote the set of neighbours of a vertex $v$ by $N(v)$ and the neighbours in $G-W$ of a vertex $v$ by $N_{G-W}(v)$. Where $S$ is a vertex set we let
\[
N(S)=\left(\cup_{v\in S}N(v)\right)\setminus S\;\text{ and }\; N_{G-W}(S)=\left(\cup_{v\in S}N_{G-W}(v)\right)\setminus S.
\] 
We will use $B(S)$ to denote the ball of radius 1 around $S$, so that $B(S)=S\cup N(S)$. Iterating this function $l$ times we get the ball of radius $l$ around $S$, denoted by $B^l(S)$. We will use similar notation like $B(v)$, $B_{G-W}(v)$ and $B^l_{G-W}(S)$ without further definition.

For a vertex set $S$, $d_S(x,y)$ is the graph theoretic distance between vertices $x$ and $y$ in the graph $G[S]$. We will denote the radius of a vertex in $S$ by $\mbox{rad}_S(v)=\max_{w\in S} d_S(v,w)$ and the radius of a set $S$ by
\[
\mbox{rad}(S)=\min\{\mbox{rad}_S(v)|v\in S\}.
\]
For a graph $G$, we will use $d(G)$ for the average degree of a graph $G$ and for a vertex set $S$ we will set $d(S)=d(G[S])$. The minimum and maximum degree of a graph $G$ will be denoted by $\d(G)$ and $\Delta(G)$ respectively.

Several lemmas require the size of a graph $H$, denoted by $m$, to be sufficiently large based on $\d$ and $t$. By this we mean there is some function $m_0(\d,t)$ for which the lemma holds if $m\geq m_0(\d,t)$.

We will use $\log$ for the natural logarithm.

\subsection{Graph Expansion}
We will use an adaptation of the expansion used by Shapira and Sudakov~\cite{SS12} to prove their result.
\de
Let $\lambda,\eta>0$. An $m$-vertex graph $H$ is said to be a \emph{$(\lambda,\eta)$-expander} if for every set $S\subseteq V(H)$ with $|S|\leq m^{1-\eta}$ we have
\[
|N(S)|\geq \lambda|S|.
\] 
\fn

We will also follow Shapira and Sudakov in using the following proposition to find an expander subgraph in graphs of constant average degree.

\begin{prop}\label{finddense}
If $G$ is a graph satisfying $d(G)=c$ and $S\subset V(G)$ satisfies $|N(S)|<\c|S|$, then either $d(G-S)\geq c$ or $d(B(S))\geq(1-\c)c$.
\end{prop}
\pr
Suppose to the contrary that $d(G-S)< c$, $d(B(S))<(1-\c)c$ and $|N(S)|<\c|S|$. Then, writing $n=|G|$,
\[
2e(G)\leq(n-|S|)c+(|S|+\c|S|)(1-\c)c=nc-\c^2|S|c<nc,
\]
a contradiction.
\oof

We will use Proposition \ref{finddense} repeatedly on graphs with constant density to find either a small dense subgraph or a larger subgraph which is an expander. In the first case the density of the subgraph will be at most a factor of $(1-\delta)$ smaller than the density of the original graph, for some fixed $\d$. In the second case if the subgraph has $m$ vertices then it will be a $(f(m),\eta)$-expander, for some fixed $\eta$. We will define an expansion function (based on parameters $\delta$ and $\eta$), so that if $f$ is such a function then this outline is possible, as shown in Lemma \ref{key}.

\de Let $G$ be a graph of order $n$ and let $0<\d,\eta<1$. A function $f:[1,n]\to[0,1]$ is a \emph{$(\d,\eta)$-expansion function for $G$} if
\begin{itemize}
\item $f(c)\leq f(a)+f(b)$ for all real $1\leq a\leq b\leq c\leq n$, and
\item $\sum_{0\leq x\leq y-2}f\left(n^{(1-\eta/2)^{x}}\right)\leq \d/2$ for $y=\lfloor-\log\log n / \log(1-\eta/2)\rfloor$.
\end{itemize}
\fn

\lem \label{key}
Let $G$ be a graph of order $n$ and let $0<\d,\eta<1$. Let $f$ be a $(\d,\eta)$-expansion function for $G$ and $c=d(G)$. Then $G$ contains a subgraph $H$ so that $d(H)\geq (1-\d)c$, $\d(H)\geq d(H)/2$ and, moreover, if $|H|> 2^{16/\eta}$, then $H$ is a $(f(|H|),\eta)$-expander.
\ma
\pr Set $G_0=G$ and consider the following process indexed by $l$. If $\d(G_l)\geq d(G_l)/2$ and if either $|G_l|\leq 2^{16/\eta}$ or $G_l$ is a $(f(|G_l|),\eta)$-expander then the process stops. Otherwise, set $m_l=|G_l|$ and $d_l=d(G_l)$. 

If $\d(G_l)<d_l/2$ then remove a vertex of minimum degree from $G_l$ to form $G_{l+1}$ with $d(G_{l+1})\geq d_l$.

If $\d(G_l)\geq d_l/2$ then $G_l$ is not a $(f(m_l),\eta)$-expander so there will be some set $S_l\subset V(G_l)$ with $|S_l|\leq m_l^{1-\eta}$ and $|N(S_l)|< f(m)|S_l|$. By Proposition \ref{finddense}, either $d(G_l-S_l)\geq d_l$ or $d(B(S_l))\geq (1-f(m_l))d_l$. In the former case let $G_{l+1}=G_l-S_l$ and in the latter case let $G_{l+1}=G[B(S_l)]$. Let $H$ be the graph at the end of this process.

We need only show that $d(H)\geq(1-\d)c$, as if $|H|> 2^{16/\eta}$ then $H$ must be an $(f(|H|),\eta)$-expander for the process to have stopped.

In the process after $l$ steps either the average degree of the graph didn't decrease, or
\[
m_{l+1}\leq (1+f(m_l))m_l^{1-\eta}\leq 2m_l^{1-\eta}< m_l^{1-\eta/2}.
\]
The density therefore can decrease only once as the size of $G_l$ passes through the interval $\left[n^{(1-\eta/2)^{x}},n^{(1-\eta/2)^{x-1}}\right]$ for each integer $x\geq 1$. The process stops once $m_l\leq e^{(1-\eta/2)^{-3}}\leq 2^{16/\eta}$, where the last inequality is true as $0<\eta<1$. We only then need to look at intervals for $x$ where $n^{(1-\eta/2)^{x-1}}\geq e^{(1-\eta/2)^{-3}}$. Taking logarithms we have that if $y=\lfloor-\log\log n / \log (1-\eta/2)\rfloor$ then $x\leq y-2$. When $1\leq x\leq y-2$ the average degree decreases as the size of $G_l$ passes through the interval for $x$ by at most a factor of $1-f(m_l)$ which is at most
\[
1-f\left(n^{(1-\eta/2)^{x}}\right)-f\left(n^{(1-\eta/2)^{x-1}}\right)
\]
by the first property in the definition of a $(\d,\eta)$-expansion function. Therefore,
\begin{align*}
d(H)&\geq \left(\prod_{x=1}^{y-2} \left(1-f\left(n^{(1-\eta/2)^{x}}\right)-f\left(n^{(1-\eta/2)^{x-1}}\right)\right)\right)c \\
&\geq \left(1-2\sum_{x=0}^{y-2}f\left(n^{(1-\eta/2)^{x}}\right)\right)c
\geq (1-\d)c,
\end{align*}
using the second property of a $(\d,\eta)$-expansion function.
\oof

The only consequence of expansion we will use can be encapsulated by the following proposition, where we expand a vertex set $S$ in $G$ while avoiding some forbidden vertex set $W$.

\begin{prop} \label{exp} Let $0<\lambda,\eta<1$. Let $H$ be a $(\lambda,\eta)$-expander of order $m$ and let $S,W\subset V(H)$ be disjoint sets with $S\neq \emptyset$. If $|W|\leq \lambda|S|/2$ and $k=(4/\lambda)\log m$ then
\[
|B_{H-W}^{k}(S)|\geq m^{1-\eta}.
\]
\end{prop}
\pr For each integer $l$, if $|B^{l}_{H-W}(S)|\leq m^{1-\eta}$, then, as $H$ is a $(\lambda,\eta)$-expander and $\lambda|B^{l}_{H-W}(S)|\geq \lambda|S|\geq 2|W|$, we have
\[
|N_{H-W}(B_{H-W}^l(S))|\geq \lambda|B_{H-W}^{l}(S)|-|W|\geq \frac12\lambda|B_{H-W}^{l}(S)|.
\]
Therefore the set $B^{l}_{H-W}(S)$ expands by at least a factor of $(1+\lambda/2)$ as $l$ increases until $l$ is the smallest integer for which $|B^{l}_{H-W}(S)|>m^{1-\eta}$. Then
\[
m\geq |B_{H-W}^{l}(S)|\geq \left(1+\frac{1}{2}\lambda\right)^{l}|S|.
\]
Taking logarithms, and using the fact that $\log (1+\lambda/2)\geq \lambda/4$ for $\lambda<1$, we get that $l\leq (4/\lambda)\log m$, as required.
\oof

The following proposition will be used several times and will be convenient to state here.

\begin{prop} \label{turan}
Suppose in a graph $G$ we have $s$ disjoint vertex sets $S_1,\ldots S_s$ with paths $P_1,\ldots,P_s$ each of length at most $t$. Then we may find a subset $I\subset \{1,\ldots,s\}$ of size at least $s/(2t+3)$ so that for every $i,j\in I$ with $i\neq j$ we have $S_i\cap P_j=\emptyset$.
\end{prop}
\pr
Form a new graph $H$ on the vertex set $\{1,\ldots,s\}$ with $ij$ an edge exactly when $S_i\cap P_j\neq\emptyset$ or $P_i\cap S_j\neq\emptyset$. Each path has at most $t+1$ vertices and thus intersects with at most $t+1$ of the disjoint sets $S_1,\ldots,S_s$. The graph $H$ thus has at most $s(t+1)$ edges. A result of Caro~\cite{caro} and Wei~\cite{wei} using a simple application of the probabilistic method says that in such a graph $\a(H)\geq \sum_{v\in V(H)}1/(d(v)+1)$. Applying Jenson's inequality, we have $\a(H)\geq s/(d(H)+1)\geq s/(2t+3)$. Picking an independent set of size at least $s/(2t+3)$ in $H$ gives the required subset $I$.
\oof

\section{Proof of Theorem \ref{logbound}}\label{main}
We will first pick many moderately-sized sets of small diameter in Proposition \ref{largesets}, before using them to construct a small $K_t$-minor in Lemma \ref{minorinexp}.
\begin{prop}\label{largesets}
The following holds for any $\lambda>0$ and constant $\eta\leq 1/8$. If $H$ is a $(\lambda,\eta)$-expander on $m$ vertices and $\lambda\geq 2m^{-1/8}$ then $H$ contains $m^{1/4}$ disjoint vertex sets $S_1,\ldots,S_{m^{1/4}}$ such that $|S_i|=m^{1/4}$ and $\mbox{rad}(S_i)\leq(4/\lambda)\log m$.
\end{prop}
\pr
Say a vertex set $S$ is \emph{nice} if $\mbox{rad}(S)\leq k:=(4/\lambda)\log m$ and $|S|=m^{1/4}$. It is enough to show that for any set $W$ of at most $m^{1/2}$ vertices, we can find in $H- W$ a nice set $S$. Indeed, we can then iteratively pick the sets $S_i$, $1\leq i\leq m^{1/4}$, where at iteration $i$ we will pick $S_i$ from $H- (\cup_{j<i}S_j)$. As $|\cup_{j<i}S_j|\leq m^{1/4}m^{1/4}=m^{1/2}$ this will be possible until we have found $m^{1/4}$ nice sets.

Given any set $W\subset V(H)$ with $|W|\leq m^{1/2}$, take a set $A\subset H-W$ with $|A|=m^{5/8}$. We have $\lambda|A|\geq2|W|$, and, hence by Proposition \ref{exp}, $|B^k_{H-W}(A)|\geq m^{1-1/8}\geq m^{1/4}|A|$. As $B^k_{H-W}(A)=\cup_{v\in A}B^{k}_{H-W}(v)$ there must be some vertex $v\in A$ for which $|B^k_{H-W}(v)|\geq m^{1/4}$. Repeatedly remove vertices from $B^k_{H-W}(v)$ of furthest distance in $H-W$ from $v$ until $m^{1/4}$ vertices remain, which will form a nice set as required.
\oof

To construct a $K_t$-minor we will begin by taking lots of moderately-sized sets $S_i$ with small diameter, provided by Proposition \ref{largesets}, and expanding them. A large number of these expanded sets will intersect, allowing us to find short paths between the sets $S_i$. We will pick a vertex $v$ in one of the original sets for which, using the small diameter of the sets, we can find short paths between $v$ and many of the sets $S_i$. We will discard the sets without a path to $v$ satisfying certain properties. Eventually, we will pick $t-1$ of the remaining paths and contract all the interior vertices into $v$ to get one vertex of the $K_t$-minor. Before we do, we put these paths to one side and repeat a similar process $t-1$ more times on the remaining sets to find other short paths which will eventually become the other vertices of the $K_t$-minor.

\lem \label{minorinexp}
The following holds for any $\lambda>0$, integer $t\geq 1$ and sufficiently large $m$ depending on $t$. If $H$ is a $(\lambda,1/8t)$-expander on $m$ vertices and $\lambda\geq 36m^{-1/8t}\log m$ then $H$ contains a $K_t$-minor consisting of at most $(16t^2/\lambda)\log m$ vertices.
\ma
\pr
Let $\eta=1/8t$. From Proposition~\ref{largesets} we know $H$ contains disjoint sets $S_1,\ldots, S_{m^{1/4}}$, each of size $m^{1/4}$ and radius at most $k:=(4/\lambda)\log m$. Pick $v_i\in S_i$ so that $d_{S_i}(v_i,y)\leq k$ for all $y\in S_i$. We will carry out the following process $t-1$ times, beginning with the indexing set $I_0=[m^{1/4}]$ and forbidden vertex set $W_0=\emptyset$.

Suppose at stage $\alpha$ we have an indexing set $I_\alpha$ with $|I_\alpha|= m^{1/4-2\alpha\eta}$ and a forbidden vertex set $W_\alpha\subset V(H)\setminus\cup_{i\in I_\alpha} S_i$ with $|W_\alpha|\leq \alpha m^{1/4-\eta}$. For sufficiently large $m$, $|W_\alpha|\leq \lambda m^{1/4}/2$. By Proposition \ref{exp} then, 
\[
|B^{k}_{H-W_\alpha}(S_i)|> m^{1-\eta}
\]
for each $i\in I_\alpha$. Therefore there must be some vertex in $H-W_\alpha$ which is in more than $|I_\alpha|m^{-\eta}$ of the sets $B^{k}_{H-W_\alpha}(S_i)$, $i\in I_\alpha$. Hence we can pick $I_\alpha'\subset I_\alpha$ so that $|I_\alpha'|\geq|I_\alpha|m^{-\eta}$ and $B^k_{H-W_\alpha}(S_i)\cap B^k_{H-W_\alpha}(S_{j})\neq\emptyset$ for each $i,j\in I_\alpha'$. 

Pick $i_\alpha\in I_\alpha'$ and remove it from $I_\alpha'$. For each remaining $i\in I_\alpha'$, let $P_{\alpha,i}$ be a path with length at most $4k$ connecting $v_{i_\alpha}$ to $v_i$ in $H-W_\alpha$ so that $P_{\alpha,i}$ is composed of two connected paths $P_{\a,i}\setminus S_{i}$ and $P_{\alpha,i}\cap S_i$ (i.e. the path does not go in and out of $S_i$). This is possible because $B^k_{H-W_\alpha}(S_i)\cap B^k_{H-W_\alpha}(S_{i_\alpha})\neq\emptyset$ for each $i\in I_\alpha'$.

Using Proposition \ref{turan}, select $I_{\alpha+1}\subset I_\alpha'$ where $|I_{\alpha+1}|=m^{1/4-2(\alpha+1)\eta}$ and if $i,j\in I_{\alpha+1}$ are distinct then $P_{\alpha,i}\cap S_{j}=\emptyset$. This is possible as, for sufficiently large $m$, $|I_\alpha'|/(8k+3)\geq m^{1/4-2(\alpha+1)\eta}$. 

Let $W_{\alpha+1}=W_{\alpha}\cup(\cup_{i\in I_{\alpha+1}}(P_{\alpha,i}\setminus S_i))$. For all distinct $i,j\in I_{\alpha+1}$, $S_{i}\cap P_{\alpha,j}=\emptyset$, so $W_{\alpha+1}$ is disjoint from each $S_i$ for $i\in I_{\alpha+1}$. Noting that, for sufficiently large $m$, $|W_{\alpha+1}|\leq |W_\alpha|+4k|I_{\alpha+1}|\leq (\alpha+1)m^{1/4-\eta}$, take $I_{\alpha+1}$ and $W_{\alpha+1}$ to be the new indexing and forbidden set respectively.

Having carried out the above process $t-1$ times the final indexing set will satisfy $|I_{t-1}|\geq m^{1/4-2(t-1)\eta}\geq 1$.

Pick any $i_t\in I_{t-1}$ and let, for each r with $1\leq r\leq t$, $X_r$ be the subgraph
\[
(\cup_{s<r}(P_{s,i_r}\cap S_{i_r})) \cup (\cup_{s>r}(P_{r,i_s}\setminus S_{i_s})).
\]
For each $r$, $1\leq r\leq t$, $X_r$ is the union of paths ending in $v_{i_r}$ and thus is a connected subgraph. For each $r$ and $s$ with $1\leq s<r\leq t$, there is an edge between $P_{s,i_r}\cap S_{i_r}\subset X_r$ and $P_{s,i_r}\setminus S_{i_r}\subset X_s$. Therefore, if for each $r$, $1\leq r\leq t$, we contract each set $X_r$ to a single vertex, using the connectivity of $X_r$, and delete all the other vertices in $H$, we get a complete subgraph on $t$ vertices. In total the subgraphs $X_r$ together have at most $4kt^2$ vertices, so we have found the $K_t$-minor we require.
\oof

The expansion function $f$ will have to fulfil three main requirements. 
Firstly, Lemma \ref{minorinexp} constructs a $K_t$-minor from $O\left(\log m/f(m)\right)$ vertices in a $(f(m),1/8t)$-expander and we wish to use it to construct a $K_t$-minor from $O(\log n)$ vertices. We thus will need $f(m)=\Omega(\log m / \log n)$.
Secondly, for any $t$, we want $f(m)\geq 36m^{-1/8t}\log m$ for sufficiently large $m$ to apply Lemma \ref{minorinexp}, which we will achieve by picking $f$ so that $f(m)=\Omega(1/(\log\log m)^2)$. Finally, we will need the function to be a $(\d,\eta)$-expansion function so that we may use Lemma \ref{key}, which we will achieve by including an appropriate factor based on the constants $\d$ and $\eta$.
To satisfy these constraints, we will take $f$ as follows.

\begin{prop}\label{expfn}
The function
\[
f(m)=\max\left\{\frac{\d\eta^2}{32(\log\log 4m)^2},\frac{\d\eta\log m}{8\log n}\right\}
\]
is a $(\d,\eta)$-expansion function. Furthermore, for sufficiently large $m$ based on $t$ we have that $f(m)\geq 36m^{-1/8t}\log m$.
\end{prop}
\pr
The first condition for $f$ to be a $(\d,\eta)$-expansion function holds as $f$ is the maximum of two monotonic functions. The final conclusion of the proposition holds as $f(m)=\Omega(1/(\log\log m)^2)$. All that remains is to check the second condition for $f$ to be a $(\d,\eta)$-expansion function. For $y=\lfloor-\log\log n / \log(1-\eta/2)\rfloor$,
\begin{align*}
\sum_{x=0}^{y-2}f\left(n^{(1-\eta/2)^{x}}\right)&\leq \sum_{x=0}^{y-2}\left(\frac{\d\eta^2}{32(\log\log n+x\log(1-\eta/2))^2}+\frac{\d\eta(1-\eta/2)^x}{8}\right)
\\
&\leq \frac{\d \eta^2}{32\log^2(1-\eta/2)}\sum_{x=0}^{y-2}\frac{1}{(y-x)^2}+\frac{\d\eta}{8}\sum^\infty_{x=0}\left(1-\frac{\eta}{2}\right)^x
\\
&\leq \frac{\d \eta^2}{32(\eta/2)^2}2+\frac{\d\eta}{8}\cdot\frac{2}{\eta}=\frac{\d}{2},
\end{align*}
where we have used that $\log (1-z)\leq -z$ for $0<z<1$.
\oof

\pr[of Theorem \ref{logbound}] Assume $\e<c(t)$, else we may take $\e=1/2$. Let $f$ be the function from Proposition \ref{expfn} with $\d=\e/3c(t)$ and $\eta=1/8t$. Let $G$ be a graph with $n$ vertices satisfying $d(G)\geq c(t)+\e$. Applying Lemma \ref{key} to G we obtain a subgraph $H$ for which
\[
d(H)\geq \left(1-\d\right)(c(t)+\e)=\left(1-\frac{\e}{3c(t)}\right)(c(t)+\e)>c(t)
\]
and, letting $m=|H|$, if $m> 2^{16/\eta}$ then $H$ is a $(f(m),\eta)$-expander. If $m> 2^{16/\eta}$ is sufficiently large for Proposition \ref{expfn} and Lemma \ref{minorinexp} to hold, say $m\geq m_0(t,\e)>2^{16/\eta}$, then we get a $K_t$-minor constructed from at most $(16t^2/f(m))\log m\leq (2^{12}t^3c(t)/\e)\log n$ vertices. If $m<m_0(t,\e)$ then using the definition of $c(t)$ we can find a $K_t$-minor in $H$, which must of course consist of at most $m_0(t,\e)$ vertices.
\oof

\section{Constants and expansion of small sets}\label{sexp2}
Constructing subdivisions is harder than constructing minors. We need to find $t$ vertices with disjoint paths joining them pairwise, rather than $t$ sets with a small radius joined up by paths, as we found when constructing minors. We will call each vertex with degree $t-1$ in a subdivision a \emph{corner vertex}. We will first find many potential corner vertices, which will have $t$ disjoint paths leading out from that vertex to $t$ disjoint sets, with these sets having medium size but a small radius. Such a structure of a vertex, $t$ paths, and $t$ sets we will call a \emph{unit}, defined explicitly later.

To find units in a graph we will need, in addition to the previous expansion properties, a minimum degree condition and a tighter control over the expansion of small sets. Roughly speaking, once we have found many units, two corner vertices can be joined up by taking a medium-sized set from each of their respective units and expanding them until the sets join up. As in the case with minors this simple idea is complicated by wishing to avoid the expansion of large sets, and so we must start with a large number of units. We will also need to be careful that the paths we are forming are disjoint.

In this section we will set out some useful constants and detail the extra conditions we need to construct a subdivision. In Section \ref{main2} we will construct units in the graph using our expansion conditions. The construction is different depending on whether a graph has lots of large degrees (Lemma \ref{manylargedegrees}) or not (Lemma \ref{fewlargedegrees}). By expanding these units we will construct subdivisions in Lemma \ref{findtop}. In Section \ref{main3} we will put this together to prove Theorems \ref{logboundtop} and \ref{topsinminors}.

\subsection{Variables and dependence}\label{constants}
We will use variables $n$, $m$, $\d$ and $\eta$ in the rest of the paper, and it will be convenient to define the variable $k$ and the function $f$ based on them here for further reference. In fact, in all cases $\eta=\frac{1}{300t^3}$, and the lemmas will eventually be applied with $\d=\frac{\e}{3}$. Note that, unlike in the previous case, $\d$ will depend only on $\e$ and not on $t$.

Given a graph $G$ with $n$ vertices we will find within it a graph $H$ with $m$ vertices which is a $(f(m),\eta)$-expander and satisfies an additional condition for the expansion of small sets detailed below. We will take $f$ to be the following function depending on $n$, $\eta$ and $\d$, which is chosen for similar reasons as the function in the minor case. Let
\begin{equation}\label{ffff}
f(m)=\max\left\{\frac{\d\eta^2}{32(\log\log 4m)^2},\frac{\d\eta\log m}{8\log n}\right\}.
\end{equation}
The variable $k$, which depends on $m$, $n$, $\d$ and $\eta$, will be the length of paths or the radius of balls. Let
\begin{equation}\label{kkkk}
k=\min\left\{\frac{128\log m(\log\log 4m)^2}{\d\eta^2},\frac{32\log n}{\d\eta}\right\}.
\end{equation}
Note that, similarly to before, $k=4\log m / f(m)$.

\subsection{Small set expansion}
\de
A graph $H$ with $m$ vertices is said to be a \emph{$(\delta,\eta,n)$-expander} if for any subset $S\subset V(H)$ with $|S|\leq m^{1-\eta}$ we have
\[
|N(S)|\geq f(m)|S|,
\]
where $f$ is the function from (\ref{ffff}) and, in addition, if $1\leq |S|\leq m^{1/3}$ then
\[
|N(S)|\geq \frac{\d}{20(\log\log 4|S|)^2}|S|.
\]
\fn

We will refer to the second condition as the expansion of small sets. We will only use it for sets of size up to $\log^2 m$, but the stronger condition will follow for no additional cost. As before, these new expander graphs can be found in graphs of constant degree.

\lem \label{key2}
Let $G$ be a graph with $n$ vertices and let $0<\d,\eta<1$. Let $f$ be the function from (\ref{ffff}) and $c=d(G)$. Then $G$ contains a subgraph $H$ so that $d(H)\geq (1-2\d)c$, $\d(H)\geq d(H)/2$ and if $|H|\geq 2^{24/\eta}$ then $H$ is a $(\d,\eta,n)$-expander.
\ma
\pr
We use the same process as in the proof of Lemma \ref{key}, the only difference being in the case where $G_{l+1}=G[B(S_l)]$. There either $|S_l|\leq m_l^{1-\eta}$ and $d(G_{l+1})\geq (1-f(m_l))d_l$, as before, or $|S_l|\leq m_l^{1/3}$ and $d(G_{l+1})\geq (1-\d/20(\log\log 4|S_l|)^2)d_l$.

The number of times the average degree can decrease and the amount by which it decreases due to the first case can be bounded as before. The second case is very similar to the first case with $\eta=2/3$. For the second case, the average degree can decrease at most once as $|G_l|$ passes through the interval $[n^{(2/3)^{z}},n^{(2/3)^{z-1}}]$ for each integer $z\geq 1$, and we only need consider intervals where $z\leq y'-2$ with $y'=\lfloor -\log\log n /\log (2/3)\rfloor$ because the process stops once $|G_l|<2^{24}<2^{24/\eta}$.

Combining the two cases we see that, with $y=\lfloor-\log\log n / \log(1-\eta/2)\rfloor$,
\begin{align*}
d(H) &\geq\left(1-2\sum_{x=0}^{y-2}f\left(n^{(1-\eta/2)^{x}}\right)-\sum_{z=0}^{y'-2}\frac{\d}{20(\log\log(n^{(2/3)^{z}})^2}\right) c 
\\
&\geq \left(1-\d-\frac{\d}{20(\log^2(2/3))}\sum_{z=0}^{y'-2}\frac{1}{(y'-z)^2}\right)c
\\
&\geq (1-2\d)c
\end{align*}
as required.
\oof

The additional condition of the expansion of small sets will be used when we have found several paths leading out of a vertex $v$ and wish to expand out from $v$ avoiding these paths.

\de
We say that paths $P_1,\ldots,P_s$, each starting with the vertex $v$ and contained in the vertex set $W$, are \emph{consecutive shortest paths from $v$ in $W$} if, for each $i$, $1\leq i\leq s$, the path $P_i$ is a shortest path between its endpoints in the set $W-\cup_{j<i}P_j+v$.
\fn

\begin{lemma} \label{corner}
Let $n,t\geq 1$, $c,\d>0$, and take $\eta=1/300t^3$. For sufficiently large $t\geq t_0(\d,c)$ the following is true. Suppose $H$ is a $(\d,\eta,n)$-expander with $m\leq n$ vertices which has a vertex $v$ with $d(v)\geq ct\sqrt{\log t}+3t$, and let $l=(\log\log m)^2$.  If $s\leq 2t$ and $P_{1},\ldots,P_{s}$ are consecutive shortest paths from $v$ in $B^{l}(v)$, and $B\subset V(H)$ with $|B|\leq t$, then
\[
|B^{l}_{H-P-B+v}(v)|\geq \log^2 m,
\]
where $P=\cup_iP_i$.
\end{lemma}
\pr Let $F=H-B-P+v$. We will show by induction on $p\geq 0$ that if $|B^p_{F}(v)|\leq \log^2 m$ then
\[
|N_{F}(B^p_{F}(v))|\geq \frac{\d}{40(\log\log4|B^p_{F}(v)|)^2}|B^p_{F}(v)|.
\]

Firstly, observe that as the paths $P_i$ are consecutive shortest paths from $v$ in $B^{l}(v)$ only the first $p+2$ vertices of $P_i$, including $v$, can belong in $N_H(B^{p}_{H-\cup_{j<i}P_j+v}(v))$. Hence
\[
|N_H(B^{p}_{F}(v))\setminus N_{F}(B^{p}_{F}(v))|\leq 2(p+1)t+t.
\]
In particular, the case $p=0$ and the degree condition for $v$ gives $|B_{F}(v)|\geq ct\sqrt{\log t}$. Let $d=ct\sqrt{\log t}$. We will consider two cases, when $|B^p_{F}(v)|$ is in the region $[d,d\log^2 d]$ and when $|B^p_{F}(v)|$ is in the region $[t\log^2 t,\log^2 m]$. As $d\geq t$, these two regions overlap to cover the region we're interested in. 

Suppose then $b=|B^p_{F}(v)| \leq d\log^2 d$. The induction hypothesis for all $0\leq p'<p$ limits how large $p$ can be. The size of the ball has increased by at most a factor of $\log^2 d$ and at each increase in radius the size increases by at least a factor of $\left(1+\d/40(\log\log (4d\log^2 d))^2\right)$ and hence
\[
p\leq \frac{\log (\log^2 d)}{\log\left(1+\d/40(\log\log (4d\log^2 d))^2\right)}
\leq \frac{200(\log\log d)^3}{\d}.
\]
Here, as previously, we have used that $\log (1+x)\geq \frac{x}{2}$ for small $x$. As $d=ct\sqrt{\log t}$, for sufficiently large $t$, we have
\[
2(p+1)t+t\leq \frac{410t(\log\log d)^3}{\d}\leq \frac{\d}{40(\log\log(4d\log^2 d))^2}d.
\]
By the expansion property then for small sets, as $b=|B^p_{F}(v)|\leq m^{1/3}$,
\begin{align*}
|N_{F}(B^p_{F}(v))|\geq \frac{\d}{20(\log\log 4b)^2}b-2(p+1)t-t\geq \frac{\d}{40(\log\log 4b)^2}b,
\end{align*}
which is the inductive hypothesis for $p$.

Now, suppose $b=|B^p_{F}(v)|\geq t\log^2 t$. The induction hypothesis for all $p'< p$ implies that
\begin{equation} \label{pupper}
p\leq \frac{\log b}{\log \left(1+\d/40(\log\log 4b)^2\right)}\leq \frac{200}{\d} (\log\log b)^2\log b.
\end{equation}
For sufficiently large $t\geq t_0(\d,c)$, we have that
\begin{equation} \label{tsl}
\frac{10^5}{\d^2}t\leq \frac{b}{(\log \log 4b)^4\log b},
\end{equation}
holds for all $b\geq t\log^2 t$. Indeed since the right hand side is increasing in $b$ we need only make $t$ sufficiently large for the inequality to hold for $b=t\log^2 t$. Then, using (\ref{pupper}) and (\ref{tsl}),
\[
2(p+1)t+t\leq \frac{410t(\log\log b)^2\log b}{\d} \leq \frac{\d}{40(\log\log 4b)^2}b. \\
\]
By the expansion property for small sets, as $b=|B^p_{F}(v)|\leq m^{1/3}$,
\begin{align*}
|N_{F}(B^p_{F}(v))|&\geq \frac{\d}{20(\log\log 4b)^2}b-2(p+1)t-t\geq \frac{\d}{40(\log\log 4b)^2}b,
\end{align*}
which is the inductive hypothesis for $p$.

The inductive hypothesis implies the set $B^p_{F}(v)$ continues to expand as $p$ increases until the first time $|B^p_{F}(v)|\geq\log^2 m$. For this value of $p$ we have $|B^{p-1}_{F}(v)|\leq\log^2 m$ so by (\ref{pupper}) we have an upper-bound for $p$. By taking sufficiently large $t$, and hence sufficiently large $m$, this bound gives $p\leq (\log\log m)^2$ as required. Here we have used that $\d$ depends only on $\e$, not on $t$.
\oof

\section{Units and Subdivisions}\label{main2}
We will first find $(t,\sigma)$-units, defined below, in expander graphs. By expanding these $(t,\sigma)$-units we can join them up to find a subdivision of $K_t$. The following definition uses the variable $k$, defined in \ref{kkkk}, and dependent on $\d$, $\eta$, $n$ and $m$. In all cases therefore we will define these variables before taking a $(t,\sigma)$-unit based on them.

\de A \emph{$(t,\sigma)$-unit} is formed from a corner vertex $v$, disjoint vertex sets $S_1,\ldots,S_{t}$ each of size $m^{\sigma}$ with special vertices $v_i\in S_i$ so that $\mbox{rad}_{S_i}(v_i)\leq k$, and a set of paths $P_1,\ldots,P_{t}$ so that $P_i$ goes from $v$ to  $v_i$, avoids all other paths $P_j$ except on $v$ itself, avoids all sets $S_j$ when $j\neq i$, and has length at most $6k$.
\fn

\subsection{Finding $(t,\sigma)$-units}
Finding $(t,\sigma)$-units requires a different technique according to whether there are many vertices of large degree (dealt with in Lemma \ref{manylargedegrees}) or few vertices of large degree (dealt with in Lemma \ref{fewlargedegrees}). In fact in this second case we may sometimes not get many $(t,\sigma)$-units, but if not we will find a small $K_t$-subdivision directly. Splitting into these two cases is governed by Proposition \ref{divide}.

\begin{prop}\label{divide}
Let $H$ be a graph with $m$ vertices and let $\sigma>0$. Then either 
\begin{enumerate}
\item there are $m^{6\sigma}$ disjoint sets $C_1,\ldots,C_{m^{6\sigma}}$ with vertices $v_i\in C_i$ such that $|C_i|=1+\log^2 m$ and $d_{C_i}(v_i)=\log^2 m$ or
\item there is a set $X$ with $|X|\leq m^{6\sigma}(1+\log^2m)$ and $\Delta(H-X)\leq \log^2 m$.
\end{enumerate}
\end{prop}
\pr If there exists $v\in V(H)$ with $d(v)\geq \log^2 m $ then let $C_1$ be $v$ along with $\log^2 m$ neighbours of $v$. Repeat this with $H-C_1$, selecting a set $C_2$ if possible, and continue. Either we will find $m^{6\sigma}$ such sets or we will find a sequence of $j$ such sets with $j\leq m^{6\sigma}$ where if $X=\cup_{i\leq j}C_i$ then $\Delta(H-X)\leq \log^2 m$.
\oof

\begin{lemma}\label{manylargedegrees}
The following holds for all $\d>0$, integers $t,n\geq 1$, and sufficiently large $m$ based on $\d$ and $t$. Let $\sigma=1/100t$ and $\eta=1/300t^3$. If $m\leq n$ and $H$ is a $(\d,\eta,n)$-expander on $m$ vertices which has $m^{6\sigma}$ sets $C_i$ containing vertices $v_i$ as in the first case of Proposition \ref{divide} then $H$ contains $m^{\sigma}$ disjoint $(t,\sigma)$-units.
\end{lemma}
\pr We will show that given any vertex set $W_0$, with $|W_0|\leq m^{3\sigma}$, we can find a $(t,\sigma)$-unit in $H-W_0$. This will prove the lemma, for if $W_0$ is taken to be the vertex set of a maximal set of disjoint $(t,\sigma)$-units in $H$ then from this $|W_0|\geq m^{3\sigma}$. Thus we must have at least $m^{\sigma}$ disjoint $(t,\sigma)$-units as required.

Given a vertex set $W_0$, with $|W_0|\leq m^{3\sigma}$, we will first expand some of the sets $C_i$ lying in $H-W_0$ to get disjoint sets $S_i$ in $H-W_0$ such that $C_i\subset S_i$, $|S_i|=m^\sigma$ and $\mbox{rad}_{S_i}(v_i)\leq k$, with $k$ defined in (\ref{kkkk}). We will say $S_i$ \emph{extends $C_i$ around $v_i$} if it has these three properties.

We will show that for any set $|W|\leq 2m^{5\sigma}$, we can find some set $C_i$ disjoint from $W$ and a set $S_i$ in $H-W$ which extends $C_i$ around $v_i$. By recursion, at each stage letting $W$ be the union of $W_0$ and the sets $S_i$ found so far, this is sufficient to find $m^{4\sigma}$ disjoint sets $S_i$ in $H-W_0$, each respectively extending $C_i$ in $H-W_0$. 

Take then $W\subset V(H)$ with $|W|\leq 2m^{5\sigma}$. Let $I$ index the sets $C_i$ which do not intersect $W$, so that $|I|\geq m^{6\sigma}-2m^{5\sigma}$. Recalling that $f=\Omega(1/(\log\log m)^2)$, we have $|\cup_{i\in I}v_i|\geq (m^{\sigma}/2-1)|W|\geq 2|W|/f(m)$ for sufficiently large $m$. By Proposition \ref{exp}, we then have
\[
|B^{k}_{H-W}(\cup_{i\in I}v_i)|\geq m^{1-\eta}\geq m^{\sigma}|\cup_{i\in I}v_i|.
\]
Thus for some $i\in I$ we have $|B^{k}_{H-W}(v_i)|\geq m^\sigma$ and so we may pick a set $S_i$ which extends $C_i$ in $H-W$ as required.

We can thus take $m^{4\sigma}$ disjoint sets $S_i$ in $H-W_0$ where each $S_i$ extends $C_i$ around $v_i$. Relabel these sets and vertices so they are indexed by $I_0=\{1,\ldots,m^{4\sigma}\}$.

Given some set of indices $I\subset I_0$, we will say that a collection of paths $\PP$ \emph{$\a$-connects the sets $S_i$, $i\in I$}, if there is a set $A=\{a_1,\ldots,a_\a\}\subset V(H)\setminus W_0$, so that $\PP$ consists of paths $P_{\beta,i}$ from $v_i$ to $a_{\beta}$ for $1\leq \beta\leq \a$ and $i\in I$, and the paths are such that
\begin{itemize}
\item $|P_{\beta,i}|\leq 2k+2$ and $|P_{\beta,i}\cap C_i|\leq 2$,
\item $P_{\beta,i}\cap W_0=\emptyset$,
\item $P_{\beta,i}\cap S_{j}=\emptyset$ if $i\neq j$, and
\item $P_{\beta,i}\cap P_{\gamma,j}=\{v_i\}\cap\{v_{j}\}$ if $\beta\neq \gamma$.
\end{itemize}

We claim that, given a collection $\PP$ that $\a$-connects the sets $S_i$, $i\in I$, where $|I|=m^{4\sigma-2\a\eta}$ and $\a<t$, we can construct a collection $\PP'$ which $(\a+1)$-connects the $S_i$ for an indexing set $I'\subset I$ with $|I'|=m^{4\sigma-2(\a+1)\eta}$.

If this claim holds we can construct a $(t,\sigma)$-unit as follows. Starting with the empty collection $0$-connecting the sets $S_i$, $i\in I_0$, and applying the claim $t$ times we get an indexing set $I_t$ and a collection $\PP$ which $t$-connects the sets $S_i$, $i\in I_t$, where $|I_t|=m^{4\sigma-2t\eta}$. As $m$ is sufficiently large we can assume without loss of generality that $\{1,\ldots,t+1\}\subset I_t$. For $1\leq s\leq t$ let $P_{s}$ be the shortest path between $v_{t+1}$ and $v_{s}$ in $P_{s,t+1}\cup P_{s,s}$. The corner vertex $v_{t+1}$, sets $S_1,\ldots, S_{t}$, and paths $P_1,\ldots,P_{t}$ form a $(t,\sigma)$-unit as required, where the properties of the $(t,\sigma)$-unit follow directly from the conditions on the paths $P_{s,i}$.

It is left then just to prove the claim. Suppose the collection $\PP$ $\a$-connects the sets $S_i$, $i\in I$, with $\a<t$. Let $W'=W_0\cup(\cup_{\beta\leq \a}\cup_{i\in I}(P_{\beta,i}-v_i))$, the vertices we wish to avoid in creating the new paths. As each path $P_{\beta,i}$ can include at most one vertex of $C_i$ in addition to $v_i$, we have $|W'\cap C_i|\leq \a\leq t$ for each $i\in I$. Therefore
\[
|\cup_{i\in I}(C_i\setminus W')|\geq |I|(\log^2m-t)\geq (3kt|I|+2|W_0|)/f(m)\geq 2|W'|/f(m),
\]
where the middle inequality follows for sufficiently large $m$ because $|W_0|\leq |I|$, $f(m)=\Omega(1/(\log\log m)^2)$ and $k=O(\log m(\log \log m)^2)$.
By Proposition \ref{exp}, 
\[
|B_{H-W'}^{k}(\cup_{i\in I}(C_i\setminus W'))|\geq m^{1-\eta}\geq m^{1/4}|\cup_{i\in I}(C_i\setminus W')|.
\]
Thus, for some $i\in I$, we have $|B_{H-W'}^{k}(C_i\setminus W'))|\geq m^{1/4}$ and hence $|B^{k+1}_{H-W'}(v_i)|\geq m^{1/4}$. Removing this $i$ from $I$ and repeating the argument, we can find a set of at least $|I|/2$ indices, indexed by $I_1\subset I$ say, so that $|B_{H-W'}^{k+1}(v_i)|\geq m^{1/4}$ for each $i\in I_1$. As $m^{1/4}f(m)\geq 2|W'|$ easily holds for large $m$, by Proposition \ref{exp} we have that, for each $i\in I_1$,
\[
|B^{2k+1}_{H-W'}(v_i)|\geq m^{1-\eta}.
\]
We can now find some vertex $a_{\a+1}$ which is in at least $|I_1|m^{-\eta}$ of these sets $B_{H-W'}^{2k+1}(v_i)$. Let the vertices $v_i$ for which this happens be indexed by $I_2$ and pick a shortest path $P_{\a+1,i}$ from $v_i$ to $a_{\a+1}$ in $H-W'$ for each $i\in I_2$. These paths will each be of length at most $2k+1$ and because each vertex in $C_i-v_i$ is adjacent to $v_i$ they will pass through $C_i-v_i$ at most once. 

Using Proposition \ref{turan}, pick an indexing set $I'\subset I_2$ of size $m^{4\sigma-2(\a+1)\eta}$ so that for each pair $i,j\in I'$ we have $S_i\cap P_{\a+1,j}=\emptyset$ if $i\neq j$. Collecting together the paths associated with each $i\in I'$ to form $\PP'$ will give the collection of paths required to finish the proof of the claim and hence the lemma.
\oof

\begin{lemma}\label{fewlargedegrees}
The following holds for all $\d,c>0$, integer $t,n\geq 1$, and sufficiently large $m$ based on $\d$, $t$ and $c$. Let $\sigma=1/100t$ and $\eta=1/300t^3$. If $m\leq n$ and $H$ is a $(\d,\eta,n)$-expander on $m$ vertices with $\d(H)\geq ct\sqrt{\log t}+3t$ and which has a set $X$ as in the second part of Proposition \ref{divide} then $H$ either contains a $K_t$-subdivision consisting of at most $t^2(\log\log m)^2$ vertices or $H$ contains $m^{\sigma}$ disjoint $(t,\sigma)$-units.
\end{lemma}

\pr Let $l=(\log \log m)^2$. We will show that for any set $W$ with $X\subset W$ and $|W|\leq m^{7\sigma}$ we can find either a $(t,\sigma)$-unit in $H-W$ or a $K_t$-subdivision with at most $t^2l$ vertices in $H$. If no such subdivision exists in $H$ then we can find $m^\sigma$ disjoint $(t,\sigma)$-units by starting with $W=X$ and repeatedly finding a $(t,\sigma)$-unit in $G-W$ to add to $W$.

Suppose then $X\subset W$ and $|W|\leq m^{7\sigma}$. Call a set $S$ \emph{nice with respect to $v\in S$} if $\mbox{rad}_S(v)\leq k$, $|S|=m^{\sigma}$ and $B_{H-W}^{l}(v)\subset S$, where $k$ is defined in \ref{kkkk}. We first wish to find a disjoint collection of $m^{1/4}$ sets in $H-W$, each of which is nice with respect to one of its vertices. If given any set $Y\subset V(H)-W$ with $|Y|\leq m^{1/2}$ we can find a nice set in $V(H)-W-Y$ then we can find $m^{1/4}$ disjoint nice sets as required. Let $Y$ then be as described.

Let $Y'=B_{H-W}^{l}(Y)$. As $\Delta(H-W)\leq\log^2 m$, 
\[
|Y'|\leq 2(\log m)^{2l}|Y|\leq 2e^{2(\log\log m)^3}|Y|.
\]
Hence, for sufficiently large $m$, $|Y'|+|W|\leq m^{3/5}$.

Take any set of vertices $V\subset V(H)-W-Y'$ with $|V|=m^{2/3}$. Recalling that $f=\Omega(1/(\log\log m)^2)$ we have, for large $m$, $2(|Y'|+|W|)\leq f(m)|V|$. By Proposition \ref{exp} then,
\[
|B^{k}_{H-Y'-W}(V)|\geq m^{\sigma}|V|,
\]
and hence there is some $v\in V$ for which $|B^{k}_{H-Y'-W}(v)|\geq m^{\sigma}$. Because $Y'=B^{l}_{H-W}(Y)$ and $v\in H-Y'-W$, we know $B^{l}_{H-W}(v)\cap Y=\emptyset$. Because $\Delta(H-W)\leq \log^2 m$, we know $|B^{l}_{H-W}(v)|\leq 2\log^{2l}m \leq m^{\sigma}$ for sufficiently large $m$. Thus we can pick $S_i\subset B^k_{H-Y-W}(v)$ so that  $B^l_{H-W}(v)\subset S_i$, $\mbox{rad}_S(v)\leq k$ and $|S|=m^{\sigma}$. Such a set $S_i$ is nice with respect to $v$ as required.

We can thus take disjoint sets $S_1,\ldots,S_{m^{1/4}}$ in $H-W$, where each set $S_i$ is nice with respect to the vertex $v_i\in S_i$.

In what follows, we will build up paths $P_{\a',i}$ and $Q_{\beta',i}$ for values of $i$, $\a'$ and $\beta'$. Paths $P_{\a',i}$ build towards constructing a $(t,\sigma)$-unit and where they cannot be found we find instead paths $Q_{\beta',i}$ which can be used directly to find a $K_t$-subdivision.

Given an indexing set $I$, we say that the collection of paths $\mathcal{P}$ \emph{$(\a,\b)$-connects the sets $S_i$}, $i\in I$, if there are sets $A=\{a_1,\ldots,a_\a\}\subset V(H)-W$ and $B=\{b_1,\ldots,b_\b\}\subset W$ so that $\mathcal{P}$ consists of paths $P_{\a',i}$ from $v_i$ to $a_{\a'}$, for $1\leq \a'\leq \a$ and $i\in I$, and paths $Q_{\b',i}$ from $v_i$ to $b_{\b'}$, for $1\leq \b'\leq \b$ and $i\in I$, and the paths are such that
\begin{itemize}
\item $|P_{\a',i}|\leq 3k$ and $|Q_{\b',i}|\leq l+1$,
\item $P_{\a',i}\subset V(H)-W$, and $Q_{\b',i}-\{b_{\b'}\}\subset B^{l}_{H-W}(v_i)$,
\item $P_{\a',i}\cap P_{\a'',j}=\emptyset$ if $i\neq j$ and $\a''\neq \a'$,
\item $P_{\a',i}\cap P_{\a'',i}=\{v_i\}$ if $\a''\neq\a'$,
\item $P_{\a',i}\cap S_j=\emptyset$ if $i\neq j$, and
\item  for each $i$ there is some ordering of the paths $P_{\a',i}$ and $Q_{\b',i}$, $1\leq \a'\leq \a$ and $1\leq \beta'\leq\beta$, so that they are consecutive shortest paths leading from $v_i$ in $B^{l}_{H-W}(v_i)$.
\end{itemize}

We claim that, given a collection $\mathcal{P}$ that $(\a,\b)$-connects the sets $S_i$, $i\in I$, $|I|\geq m^{1/4-8(\a+\b)\eta}$ with $\a,\b\leq t$ we can construct a collection $\mathcal{P}'$ and an indexing set $I'$ with $|I'|\geq m^{1/4-8(\a+\b+1)\eta}$ so that either $\mathcal{P}'$ $(\a+1,\b)$-connects the sets $S_i$, $i\in I'$, or $\mathcal{P}'$ $(\a,\b+1)$-connects the sets $S_i$, $i\in I'$.

To prove the claim, first let $P$ be the union of all the paths in $\mathcal{P}$ and note that for each $i\in I$ we have either
\begin{enumerate}
\item there is a path of length at most $l$ from $v_i$ which ends in $W-B$ and otherwise lies entirely within $B^{l}_{H-W-P+v_i}(v_i)\subset S_i$, 
\item or $B_{H-W-P+v_i}^{l}(v_i)=B^{l}_{H-B-P+v_i}(v_i)=B^{l}_{H-B-P\cap S_i+v_i}(v_i)$. In this second case by Lemma \ref{corner} applied to the paths $P\cap S_i$ we have then
\[
|B_{H-W-P+v_i}^{l}(v_i)|\geq \log^2m.
\]
\end{enumerate}

If we have the first case for at least $m^{1/4-8(\a+\b)\eta}/2$ values of $i\in I$, then there must be some vertex in $W$, $b_{\b+1}$ say, which is the end of at least $m^{1/4-8(\a+\b)\eta-7\eta}/2$ of these paths, indexed say by $I'$. Denote the path from $v_i$ to $b_{\b+1}$ by $Q_{\b+1,i}$ for each $i\in I'$ and add them to the other paths indexed by $i\in I'$ to form $\mathcal{P}'$ which $(\a,\b+1)$-connects the sets $S_i$, $i\in I$.

Suppose then we have the second case for at least $m^{1/4-8(\a+\b)\eta}/2$ values of $i\in I$, say indexed by $I_1$. Take $T_i\subset B_{H-W-P+v_i}^{l}(v_i)$ so that $v_i\in T_i$, $\mbox{rad}_{T_i}(v_i)\leq l$ and $|T_i|=\log^2 m$. Remove from $P$ the vertices $v_i$, $i\in I_1$, and the vertices which are only in paths indexed by $I\setminus I_1$. We have $|\cup_{i\in I_1}T_i|\geq \log^2 m |I_1|$, $|P|\leq 4tk|I_1|$ and $|W|\leq m^{7\sigma}$. As $k=O(\log m(\log\log m)^2)$ we have, for sufficiently large $m$, that $f(m)|\cup_{i\in I_1}T_i|\geq 2|W\cup P|$. By Proposition \ref{exp} then,
\[
|B^{k}_{H-W-P}(\cup_{i\in I_1}T_i)|\geq m^{1/4}|\cup_{i\in I_1}T_i|.
\]
There must therefore be some $i\in I_1$ for which $|B^{k}_{H-W-P}(T_i)|\geq m^{1/4}$. Remove $i$ from $I'$. Doing this repeatedly allows us to find a set $I_2\subset I_1$ with $|I_2|\geq m^{1/4-8(\a+\b)\eta}/4$ such that $|B^{k}_{H-W-P}(T_i)|\geq m^{1/4}$ for all $i\in I_2$. Applying Proposition \ref{exp} again shows that $|B^{2k}_{H-W-P}(T_i)|\geq m^{1-\eta}$ for all $i\in I_2$. There must then be some vertex $a_{\a+1}$ in $H-W-P$ which is in at least $m^{1/4-8(\a+\b)\eta-\eta}/4$ of the sets $B^{2k}_{H-W-P}(T_i)$. Let $P_{\a+1,i}$ be a shortest path in $H-W-P$ from $v_i$ to $a_{\a+1}$, which will have length at most $2k+l$. Using Proposition \ref{turan}, select $m^{1/4-8(\a+\b+1)\eta}$ of these indices, indexed by $I'$, so that $P_{\a+1,i}$ does not intersect $S_j$ if $i\neq j$. Removing paths associated with indices $I\setminus I'$ from $\mathcal{P}$ and adding in the paths $P_{\a+1,i}$ for $i\in I'$ gives $\mathcal{P}'$ which $(\a+1,\b)$-connects the sets $S_i$, $i\in I$. This completes the proof of the claim.

Starting then with $I=[m^{1/4}]$ and $\mathcal{P}=\emptyset$ which $(0,0)$-connects the sets $S_i$, $i\in I$, apply the claim at most $2t$ times to get $\mathcal{P}$ and $I'$ so that $\mathcal{P}$ $(\a,\b)$-connects the sets $S_i$, $i\in I'$, with either $\a= t$ and $\b\leq t$ or $\a\leq t$ and $\b=t$, and $|I|'\geq m^{1/4-8(\a+\b)\eta}\geq m^{1/4-16t\eta}$. For sufficiently large $m$, $|I'|\geq \binom t 2$. Without loss of generality, suppose $\{1,\ldots,\binom t 2\}\subset I'$.

If $\a=t$, then taking $v_{t+1}$ as the corner vertex, and for each $i$, $1\leq i\leq t$, the set $S_i$ and a path $P_i$ of length at most $6k$ connecting $v_i$ to $v_{t+1}$ through $P_{i,i}$ and $P_{i,t+1}$, the properties of the paths $P_{i,j}$ listed above imply that this is a $(t,\sigma)$-unit.

If $\b= t$, then let $c(i,j):[t]^{(2)}\to \{1,\ldots,\binom t2\}$ be a bijection. When $1\leq i<j\leq t$, let $R_{i,j}$ be a path of length at most $2l$ from $b_i$ to $b_j$ in $Q_{i,c(i,j)}\cup Q_{j,c(i,j)}$. From the properties of the paths $Q_{i,j}$ we have that the paths $R_{i,j}$ intersect only at their endpoints and hence the paths $R_{i,j}$ form a $K_t$-subdivision in $H$ with corner vertices $\{b_1,\ldots,b_t\}$ which contains at most $2\binom{t}{2}l+t\leq t^2l$ vertices in total, as required.
\oof

\subsection{Constructing subdivisions from $(t,\sigma)$-units}
\begin{lemma} \label{findtop}
The following holds for all $\d>0$, integer $t,n\geq 1$, and sufficiently large $m$ based on $\d$ and $t$. Let $\sigma=1/100t$ and $\eta=1/300t^3$. If $m\leq n$ and $H$ is a $(\d,\eta,n)$-expander of order $m$ which has $m^{\sigma}$ disjoint $(t,\sigma)$-units then it contains a $K_t$-subdivision consisting of at most $(10^6t^5/\d)\log n$ vertices.
\end{lemma}
\pr
Label the $(t,\sigma)$-units $U_i$, $1\leq i\leq m^{\sigma}$, where each unit $U_i$ has a corner vertex $v_i$ and paths $P_{i,j}$ from $v_{i}$ to the vertex $v_{i,j}$ in the set $S_{i,j}$, for $1\leq j\leq t$. Delete vertices from each set $S_{i,j}$ so that $|S_{i,j}|=m^{\sigma-(j-1)\eta}$ yet still $\mbox{rad}_{S_{i,j}}(v_{i,j})\leq k$.

Consider the lexiographic ordering on ordered pairs $(\a,\b)\subset [t]^2$, so that $(\a',\b')\leq (\a,\b)$ if either $\a'<\a$ or $\a'=\a$ and $\b'\leq \b$. We will say $(\a,\b)^+=(\a^+,\b^+)$ is the successor of $(\a,\b)$ under this ordering of $[t]^2$. We will also say that $(1,0)\leq (\a,\b)$ for all $(\a,\b)\subset [t]^2$ with $(1,0)^+=(1,1)$.

Given an indexing set $I$, we say that a collection of paths $\QQ$ \emph{$(\a,\b)$-connects the units $U_i$}, $i\in I$, if it consists of paths $Q_{\a',\b',i}$ for $(1,1)\leq (\a',\b')\leq (\a,\b)$ and $i\in I$ and there is a set of vertices $\{a_{\a',\b'}:(1,1)\leq (\a',\b')\leq (\a,\b)\}$ such that
\begin{itemize}
\item $Q_{\a',\b',i}$ connects $v_{i,\a'}$ to $a_{\a',\b'}$,
\item $Q_{\a',\b',i}$ has length at most $2k$,
\item $Q_{\a',\b',i} \cap Q_{\a'',\b'',j}=\emptyset$ if $\a'\neq \a'$, or if $\b'\neq\b''$ and $i\neq j$,
\item $Q_{\a',\b',i}\cap P_{j,\a''}=\emptyset$ unless $i=j$ and $\a'=\a''$, and
\item $Q_{\a',\b',i} \cap S_{j,\a''}=\emptyset$ if $\a''> \a'$ or if $\a''=\a'$ and $i\neq j$.
\end{itemize}
Note that this is a different definition to the connection of vertex sets in the proof of Lemma \ref{fewlargedegrees}. Roughly speaking, the first four conditions give paths of length at most $2k$ which connect $v_{i,\a'}$ to $a_{\a',\b'}$ while avoiding all other paths except for those also coming from $v_{i,\a'}$ or also going to $a_{\a',\b'}$. Thus each vertex $a_{\a',\b'}$ will give us one opportunity to create a path from some vertex $v_{i,\a'}$ to some other vertex $v_{j,\a'}$ which, when combined with $P_{i,\a'}$ and $P_{j,\a'}$, will allow us to create a path from $v_i$ to $v_j$. The final condition dictates the paths avoid enough of the sets $S_{j,\a''}$ to prove the following claim.

We claim that, given a collection $\QQ$ that $(\a,\b)$-connects the $(t,\sigma)$-units $U_i$, $i\in I\subset\{1,\ldots,m^\sigma\}$, where $|I|= m^{\sigma-2\eta(\a t+\b)}$, we can construct a collection $\QQ'$ and indexing set $I'$ with $|I'|=m^{\sigma-2\eta(\a^+ t+\b^+)}$ so that $\QQ'$ $(\a,\b)^+$-connects the units $U_i$, $i\in I'$.

Given this claim, we can construct a subdivision on $O_{\e,t}(\log n)$ vertices as follows. Discard some of the units $U_i$ so that $m^{\sigma-2\eta t}$ $(t,\sigma)$-units remain, indexed by $I$. Starting with the empty collection which $(1,0)$-connects the units $U_i$, $i\in I$, apply the claim $t^2$ times to get a collection $\QQ$ and an indexing set $I'$ so that $\QQ$ $(t,t)$-connects the units $U_i$, $i\in I'$, and $|I'|=m^{\sigma-2(t^2+t)\eta}$. Taking $m$ to be sufficiently large, we can assume without loss of generality that $\{1,\ldots,t\}\subset I'$.

It is well known that $K_t$ can be $t$ edge-coloured with colours $1,\ldots,t$ so that edges of the same colour do not meet. Take such a colouring $c:E(K_t)\to[t]$. There are certainly at most $t$ edges of each colour, so we may construct a function $e:E(K_t)\to [t]$ such that if $ij\neq i'j'$ then $c(ij)\neq c(i'j')$ or $e(ij)\neq e(i'j')$. For $i<j$, let $R_{ij}$ be a path of length at most $16k$ connecting $v_i$ to $v_j$ in $P_{i,c(ij)}\cup Q_{c(ij),e(ij),i}\cup Q_{c(ij),e(ij),j}\cup P_{j,c(ij)}$. By the properties of these paths and the functions $c$ and $e$ we know that $R_{ij}\cap R_{i'j'}=\{v_i,v_j\}\cap\{v_{i'},v_{j'}\}$. Therefore the paths $R_{i,j}$ form a $K_t$-subdivision with the vertices $v_i$ as corner vertices which uses at most $16kt^2\leq (10^6t^5/\d)\log n$ vertices as required.

It is left then to prove the claim. Suppose $\QQ$ $(\a,\b)$-connects the $(t,\sigma)$-units $U_i$, $i\in I$, where $|I|= m^{\sigma-2\eta(\a t+\b)}\leq m^{\sigma-2\eta\a}$. Let $W$ be the union of all the paths $Q_{\a',\b',i}$ in $\QQ$ and the paths $P_{i,\a'}$ with $1\leq \a'\leq t$ and $i\in I$, so that, for sufficiently large $m$, $|W|\leq m^{\sigma-2\eta\a}(t^2+t)(6k+1)\leq f(m)m^{\sigma-\a\eta}/4$. Let $W_i=(W\setminus S_{i,\a^+})\cup(\cup_{\a'>\a^+}S_{i,\a'})$, so that
\[
|W_i|\leq |W|+\sum_{\a'>\a^+}m^{(\sigma-\a')\eta} \leq |W| + tm^{\sigma-(\a^++1)\eta}\leq \frac{f(m)}{2}m^{\sigma-\a^+\eta}.
\]
Then, by Proposition \ref{exp}, we have for each $i\in I$
\[
|B^k_{H-W_i}(S_{i,\a^+})|\geq m^{1-\eta}.
\]
Therefore, there must be some vertex, $a_{\a^+,\b^+}$ say, which is in at least $m^{-\eta}|I|$ of the sets $B^k_{H-W_i}(S_{i,\a^+})$, indexed say by $I_1$. Pick paths $Q_{\a^+,\b^+,i}$ of length at most $2k$ from $v_{i,\a}$ to $a_{(\a,\b)^+}$ in $H-W_i$ for each $i\in I_1$. Using Proposition \ref{turan} on the paths $Q_{\a^+,\b^+,i}$ and matching sets $\cup_{\a\geq \a^+}S_{i,\a}$ we can find a further indexing set $I'\subset I_1$ so that for all $i,j\in I'$ we have $(\cup_{\a\geq \a^+}S_{i,\a})\cap Q_{\a^+,\b^+,j}=\emptyset$ if $i\neq j$ and $|I'|= m^{\sigma-2\eta(\a t+\b+1)}$. As $\a t+\b+1=\a^+t+\b^+$, this completes the proof of the claim, and hence the lemma.
\oof

\section{Proofs of Theorems \ref{logboundtop} and \ref{topsinminors}}\label{main3}
Theorems \ref{logboundtop} and \ref{topsinminors} follow from Lemma \ref{key2} and the following lemma.

\lem \label{minorinexp2}
Let $\d,c>0$. For sufficiently large $t\geq t_0(\d,c)$, $\eta=1/300t^3$, and sufficiently large $m\geq m_0(\d,t,c)$ the following holds. If $H$ is a $(\d,\eta,n)$-expander on $m$ vertices, with $\d(H)\geq ct\sqrt{\log t}+3t$, then $H$ has a $K_t$-subdivision consisting of at most $(10^6t^5/\d)\log n$ vertices.
\ma

\pr Lemma \ref{minorinexp2} follows simply from Proposition \ref{divide} and Lemmas \ref{manylargedegrees}, \ref{fewlargedegrees} and \ref{findtop}.
\oof

\pr[of Theorem \ref{logboundtop}] Assume $\e<1/2$, for otherwise we may set $\e=1/3$. Let $\d=\e/3$, so that $(1-2\d)(1+\e)>1$. Take $t\geq t_0(\d,1/128)$ where $t_0$ comes from Lemma \ref{minorinexp2}.

Let $G$ be a graph with $n$ vertices satisfying $d(G)\geq s(t)(1+\e)$. Applying Lemma \ref{key2} to $G$ we obtain a subgraph $H$ so that
\[
d(H)\geq (1-2\d)(1+\e)s(t)>s(t),
\]
$\d(H)\geq s(t)/2\geq t^2/128$, and if $|H|> 2^{24/\eta}$ then $H$ is a $(\d,\eta,n)$-expander, where $\eta=1/300t^3$. From Lemma \ref{minorinexp2}, if $|H|\leq\max\{m_0(\d,t,1/128),2^{24/\eta}\}$ then from the definition of $s(t)$ we can find a constant-sized $K_t$-subdivision, and otherwise we can construct a $K_t$-subdivision from $O_{\e,t}(\log n)$ vertices using Lemma \ref{minorinexp2}.
\oof

\pr[of Theorem \ref{topsinminors}] Assume $\e<1/2$, for otherwise we may set $\e=1/3$. Let $\d=\e/3$, so that $(1-2\d)(1+\e)>1$. Take $t\geq t_0(\d,1/128)$ where $t_0$ comes from Lemma \ref{minorinexp2}.

Let $G$ be a graph with $n$ vertices satisfying $d(G)\geq c(t)(1+\e)$. Applying Lemma \ref{key2} to G we obtain a subgraph $H$ satisfying
\[
d(H)\geq (1-2\d)(1+\e)c(t)>c(t),
\]
$\d(H)\geq c(t)/2\geq t\sqrt{\log t}/128+3t$ and if $|H|>2^{24/\eta}$ then $H$ is a $(\d,\eta,n)$-expander, where $\eta=1/300t^3$. From Lemma \ref{minorinexp2}, if $|H| \leq\max\{m_0(\d,t,1/128),2^{24/\eta}\}$, then from the definition of $c(t)$ we can find a constant-sized $K_t$-minor, and otherwise we can construct a $K_t$-subdivision from $O_{\e,t}(\log n)$ vertices using Lemma \ref{minorinexp2}.
\oof

\section{Final remarks}\label{remarks}
For the proof of Theorems \ref{logboundtop} and \ref{topsinminors} we used that the extremal functions for topological minors and minors were of order $\Theta(t\sqrt{\log t})$ and $\Theta(t^2)$ respectively. Specifically, they were used in the construction of the corner vertices in the $(t,\sigma)$-units which had $t$ disjoint paths emerging from them. We were able to construct these paths and hence the units as the minimum degree was comfortably larger than just a linear function in $t$. If we could have constructed $\binom t 2$ disjoint paths leading out of the corner vertices to larger sets the proof of Lemma \ref{findtop} would be more straightforward but the degree conditions are not strong enough to allow that.

An \emph{immersion minor} is a subdivision where in addition we allow the paths between the corner vertices to intersect on interior vertices, only requiring them to be edge-disjoint. The extremal function for immersion minors was shown to be $\Theta(t)$ by DeVos, Dvo{\v{r}}{\'a}k, Fox, McDonald, Mohar and Scheide~\cite{immersion}. This is not large enough to apply our methods to find a corresponding version of Theorem \ref{logbound} or Theorem \ref{logboundtop} for immersion minors.

\subsection*{Acknowledgements}
The author would like to thank Andrew Thomason for his help and suggestions.

\bibliographystyle{plain}
\bibliography{rhmreferences}

\end{document}